\documentclass[10pt]{article}
\usepackage{latexsym}
\usepackage{amsfonts}
\usepackage{enumerate}
\usepackage{multicol}
\usepackage{graphicx}
\usepackage{amssymb}
\usepackage{amsmath}
\usepackage{epic}

\title{On the Diameter of Undirected Cayley Graphs of Finite Abelian Groups\\[.4in]}

\author{B\'{e}la Bajnok\footnote{Department of Mathematics, Gettysburg College, U.S.A.  Email: bbajnok@gettysburg.edu.} 
\hspace{.7in} W. Kyle Beatty\footnote{Department of Applied Mathematics and Statistics, Johns Hopkins University, U.S.A.  Email: wbeatty1@jhu.edu.}}

\date{April 1, 2024}

\newtheorem{thm}{Theorem}[section]

\newtheorem{lem}[thm]{Lemma}
\newtheorem{cor}[thm]{Corollary}
\newtheorem{prop}[thm]{Proposition}
\newtheorem{conj}[thm]{Conjecture}

\begin{document}

\maketitle

\begin{abstract}  Let $s$ be a positive integer.  Our goal is to find all finite abelian groups $G$ that contain a $2$-subset $A$ for which the undirected Cayley graph $\Gamma(G,A)$ has diameter at most $s$.  We provide a complete answer when $G$ is cyclic, and a conjecture and some partial answers when $G$ is noncyclic.  

\medskip


\medskip

      
\end{abstract}

\thispagestyle{empty}

\section{Introduction}

Throughout this paper we let $G$ be an additively written finite abelian group.  When $G$ is cyclic and of order $n$, we identify it with $\mathbb{Z}_n=\mathbb{Z}/n\mathbb{Z}$; we consider $0,1,\dots,n-1$ interchangeably as integers and as elements of $\mathbb{Z}_n$.  Below we use standard terminology and notation of additive combinatorics (see \cite{Baj:2018a}).

Given a finite abelian group $G$ and a subset $A \subseteq G$, the {\em undirected Cayley graph} $\Gamma (G,A)$ has its vertex set formed by the elements of $G$, and has an edge between two vertices whenever their difference is an element of $A$, that is, $g_1 \in G$ and $g_2 \in G$ are connected when $g_2=g_1+a$ or $g_1=g_2+a$ for some $a \in A$.  If $A$ is a generating set of $G$, then $\Gamma (G,A)$ is connected; in particular, if $A$ is an {\em $s$-spanning set} of $G$, meaning that every element of $G$ can be written as the signed sum of at most $s$ (not-necessarily distinct) elements of $A$, then $\Gamma(G,A)$ has diameter at most $s$.  

Suppose now that $A$ is an $m$-subset of $G$.  It is easy to compute that there are at most $\sum_{i \geq 0} {m \choose i} {h-1 \choose i-1}2^i $ vertices of $\Gamma(G,A)$ that are of distance $h$ away from any given vertex, and thus if $A$ is an $s$-spanning set in $G$, then the order of $G$ can be at most
$$\sum_{h=0}^s \sum_{i \geq 0}  {m \choose i} {h-1 \choose i-1} 2^i = \sum_{i \geq 0}  {m \choose i} {s \choose i}2^i .$$ 
Therefore, we have the following bound.

\begin{prop} \label{diameter a(m,s)}
For positive integers $m$ and $s$, define 
\begin{eqnarray} \label{eqn for a(m,s)} a(m,s)= \sum_{i \geq 0}  {m \choose i} {s \choose i}2^i .\end{eqnarray}
If an abelian group $G$ has an $m$-subset $A$ for which the diameter of $\Gamma (G,A)$ is $s$, then the order of $G$ is at most $a(m,s)$.

\end{prop}
We note, in passing, that the values of $a(m,s)$ are known as {\em Delannoy numbers} (see \cite{Del:1895a}), and that they satisfy the recursion
$$a(m,s)=a(m,s-1)+a(m-1,s)+a(m-1,s-1).$$ 

It is a natural question to ask for the largest possible order of a group $G$ with an $m$-subset $A$ for which the diameter of  ${\Gamma} (G,A) $ is at most $s$.      Before we discuss our results on this question, let us review the analogous question in the directed case.  

Given a finite abelian group $G$ and a subset $A \subseteq G$, the {\em directed Cayley graph} $\vec{\Gamma} (G,A)$ has its vertex set formed by the elements of $G$, and has a directed edge from $g_1 \in G$ to $g_2 \in G$ when $g_2=g_1+a$ for some $a \in A$.  We say that $A$ is an {\em $s$-basis} of $G$ if every element of $G$ can be written as the sum of at most $s$ (not-necessarily distinct) elements of $A$; in this case, the diameter of $\vec{\Gamma} (G,A)$ is at most $s$.

For an $m$-subset $A$ of $G$, we see that there are at most ${m+h-1 \choose h}$ vertices of $\vec{\Gamma}(G,A)$ that are of distance $h$ away from any given vertex, so if $A$ is an $s$-basis of $G$, then we must have 
$$|G| \leq \sum_{h=0}^s {m+h-1 \choose h} = {m+s \choose s}.$$

 When equality holds, we say that  $A$ is a {\em perfect $s$-basis} of $G$: in this case, every element of $G$ arises as the sum of at most $s$ (not-necessarily distinct) elements of $A$ in a unique way.  Trivially, the set of nonzero elements is a perfect $1$-basis in any
group $G$, and the $1$-element set consisting of a generator of $G$ is a perfect $s$-basis
in the cyclic group of order $s +1$. It turns out that there are no others:

\begin{thm} [Bajnok, Berson, and Just; cf. \cite{BajBerJus:2022a}] \label{BajBerJus:2022a} If a subset $A$ of a finite abelian group $G$ is a perfect $s$-basis, then $s =1$
and $A = G \setminus \{0\}$, or $G \cong \mathbb{Z}_{s+1}$ and $|A| = 1$.

\end{thm}

Given the scarcity of perfect bases, it is natural to ask for the largest groups with a given diameter and size of generating set.  The paper \cite{MasSchJia:2011a} asks the following questions: 
\begin{itemize}
  \item What is the largest possible order of an abelian group $G$ with an $m$-subset $A$ for which the diameter of  $\vec{\Gamma} (G,A) $ is at most $s$?
  \item What is the largest possible order of a cyclic group $G$ with an $m$-subset $A$ for which the diameter of $\vec{\Gamma} (G,A)$ is at most $s$?
\end{itemize}

After a series of papers \cite{MorFioFab:1985a, FioYebAleVal:1987a, HsuJia:1994a, MasSchJia:2011a, Fio:2013a}, these two questions have been completely solved for the case of $m=2$.  We have the following two results.

\begin{thm} [Morillo, Fiol, and F\`abrega; cf.~\cite{MorFioFab:1985a}] \label{thmnform=2} 

Let $s$ be a positive integer.  The largest possible order of a cyclic group $G$ with a $2$-subset for which the diameter of $\vec{\Gamma} (G,A)$ is at most $s$ is given by $\left\lfloor (s^2+4s+3)/3 \right\rfloor$.

\end{thm}

\begin{thm} [Fiol; cf.~\cite{Fio:2013a}] \label{Fiol 2013}
Let $s$ be any positive integer.  The largest possible order of an abelian group $G$ with a $2$-subset for which the diameter of  $\vec{\Gamma} (G,A) $ is at most $s$ is given by $\left\lfloor (s^2+4s+4)/3 \right\rfloor$.

\end{thm}

Observe that the quantities in Theorems \ref{thmnform=2} and \ref{Fiol 2013} agree, unless $s \equiv 1$ mod 3; for that case one can verify that the group $G=\mathbb{Z}_k \times \mathbb{Z}_{3k}$ has the $s$-basis $\{(0,1),(1,3k-1)\}$, where $k=(s+2)/3$.  However, the question specifically for noncyclic groups is still open when $s \not \equiv 1$ mod 3.

Let us now return to the undirected case.  When equality occurs in Proposition \ref{diameter a(m,s)}, then $A$ is called a {\em perfect $s$-spanning set} in $G$.  Perfect spanning sets can be easily characterized for $m=1$ and $s=1$.  
Note that $a(1,s)=2s+1$ and, clearly, for any $a \in \mathbb{Z}_{2s+1}$ with $\gcd(a,2s+1)=1$, the set $\{a\}$ is a perfect $s$-spanning set in $\mathbb{Z}_{2s+1}$.  Also, $a(m,1)=2m+1$, and for any group $G$ of order $2m+1$, every asymmetric subset of size $m$ (that is, a set containing exactly one of any nonzero element or its inverse) is a perfect $1$-spanning set in $G$.  

From (\ref{eqn for a(m,s)}) we get $a(2,s)=2s^2+2s+1$, and below we establish the following result.

\begin{thm} \label{perfect span for m=2} 
Let $s$ be a positive integer.  The largest possible order of a cyclic group $G$ with a $2$-subset for which the diameter of ${\Gamma} (G,A)$ is at most $s$ is given by $2s^2+2s+1$.

\end{thm}

For noncyclic groups the situation seems a bit more nuanced.  We have the following conjecture.

\begin{conj} \label{conj on f pm (2,s) noncyclic}
Let $s \geq 3$ be an integer.  
\begin{itemize}
  \item If $2s+1$ is prime, then the largest noncyclic group with an $s$-spanning set of size $2$ has order
$2s^2$.  
\item If $2s+1$ is composite and $p$ is its smallest prime divisor, then the largest noncyclic group with an $s$-spanning set of size $2$ has order
$2s^2 +2s - (p^2-1)/2$.
\end{itemize}
In particular, for any noncyclic group $G$ with an $s$-spanning set of size $2$ we have $|G| \leq 2s^2+2s-4.$

\end{conj}
(We add that the largest noncyclic group with a $2$-spanning set of size $2$ is $\mathbb{Z}_3 \times \mathbb{Z}_3$.)

Below we prove that the quantities in Conjecture \ref{conj on f pm (2,s) noncyclic} are lower bounds. 

\begin{thm}
Let $s$ be a positive integer.  
\begin{itemize}
  \item The largest noncyclic group with an $s$-spanning set of size $2$ has order at least $2s^2$.  
  \item If $2s+1$ is composite and $p$ is any prime divisor of $2s+1$, then the largest noncyclic group with an $s$-spanning set of size $2$ has order at least $2s^2 +2s - (p^2-1)/2$.
\end{itemize}

\end{thm}

Establishing upper bounds appears to be more difficult, and we only have the partial results that we summarize here.  We first note that a noncyclic group that has a spanning set of size $2$ must have rank $2$ and thus is isomorphic to $\mathbb{Z}_c \times \mathbb{Z}_{kc}$ for some integers $c \geq 2$ and $k \geq 1$.  

We say that the group $\mathbb{Z}_c \times \mathbb{Z}_{kc}$ is {\em $s$-regular} for some positive integer $s$, if it possesses an $s$-spanning set of the form $\{(1,u),(1,v)\}$ for some $u, v \in \mathbb{Z}_{kc}$.  While we have no precise count on the number of $s$-regular groups, it appears that most groups of rank $2$ that have an $s$-spanning set of size $2$ are $s$-regular.  For example, for $s=10$, of the $105$ groups of rank $2$ that have an $s$-spanning set of size $2$, $103$ are $s$-regular.  (The two exceptions are $\mathbb{Z}_8 \times \mathbb{Z}_{16}$ and $\mathbb{Z}_{10} \times \mathbb{Z}_{20}$.) In Theorem \ref{s regular prime c} below we prove for every odd prime $c$ and positive integer $k$ that if $\mathbb{Z}_c \times \mathbb{Z}_{kc}$ has an $s$-spanning set at all, then it is  $s$-regular.

Our results for upper bounds can be summarized as follows.   

\begin{thm}
Let $s$ be a positive integer, and let $G=\mathbb{Z}_c \times \mathbb{Z}_{kc}$ be an $s$-regular group for some integers $c \geq 2$ and $k \geq 1$.  
\begin{itemize}
  \item If $c=3$ and $s+1$ is divisible by $3$, then $|G| \leq 2s^2+1$.
  \item If $c$ is odd and $s$ is divisible by $c$, then $|G| \leq 2s^2+s$, and if $2s+1$ is divisible by $c$, then $|G| \leq  2s^2+2s-(c^2-1)/2$.
  \item In all other cases, $|G| \leq  2s^2$.
\end{itemize}

\end{thm}
As a consequence of our results, we can prove that for $s \geq 3$, every $s$-regular abelian group of rank 2 has order at most $2s^2+2s-4$, with equality if, and only if, it is isomorphic to $\mathbb{Z}_3 \times \mathbb{Z}_{(2s^2+2s-4)/3}$ for some $s \equiv 1$ mod $3$ (See Theorem \ref{thm on s-regular max order} below).  

For $c=2$ we have exact results, even without the assumption that the group is $s$-regular:

\begin{thm}
Let $k$ and $s$ be positive integers.  The largest possible order of a group of the form $\mathbb{Z}_2 \times \mathbb{Z}_{2k}$ that has an $s$-spanning set of size $2$ is $2s^2$ when $s$ is even and $2s^2-2$ when $s$ is odd.

\end{thm} 
In Section \ref{section groups 2 by 2k} below we provide some further results on groups of the form $\mathbb{Z}_2 \times \mathbb{Z}_{2k}$.

For $m \geq 3$, the largest order of abelian groups with an $m$-subset for which the diameter of their Cayley graphs is at most $s$ is generally unknown, though some small cases for cyclic groups were displayed by Haanp\"a\"a in \cite{Haa:2004a}, by Hsu and Jia in \cite{HsuJia:1994a}, and by Graham and Sloane in \cite{GraSlo:1980a}.  

Of particular interest are perfect spanning sets.  An  $m$-subset $A$ of $G$ is called a {\em perfect $s$-spanning set}, if every element of $G$ arises as a signed sum of at most $s$ terms in a unique way; in this case we have $|G|=a(m,s)$.  Above we mentioned perfect spanning sets for $m=1$, for $s=1$, and for $m=2$; we believe that there are no others.

\begin{conj}
If an $m$-subset $A$ of a finite abelian group $G$ is a perfect $s$-spanning set, then one of the following holds:
\begin{itemize}
  \item $m =1$, $G \cong \mathbb{Z}_{2s+1}$, and $A = \{a\}$ for some $a \in \mathbb{Z}_{2s+1}$ with $\gcd(a,2s+1)=1$; 
  \item $s=1$, and $G$ is the disjoint union of $\{0\}$, $A$, and $-A$; or
  \item $m=2$, $G \cong \mathbb{Z}_{2s^2+2s+1}$, and $A=d \cdot \{s,s+1\}$ for some $d \in \mathbb{Z}_{2s+1}$ with $\gcd(d,2s^2+2s+1)=1$.
\end{itemize}

\end{conj}

In Sections \ref{section Constructions and lower bounds} and \ref{section -regular groups and upper bounds} below, we present our results on lower bounds and upper bounds, respectively. In Section \ref{section groups 2 by 2k} we treat the case of the groups $\mathbb{Z}_2 \times \mathbb{Z}_{2k}$ in more detail.

\section{Constructions and lower bounds}  \label{section Constructions and lower bounds}

In this section we provide some explicit constructions for $s$-spanning sets of size $2$.  We denote the $s$-span of a 2-subset $A=\{a_1,a_2\}$ of $G$ by
$$\langle A \rangle_s = \{ \lambda_1 a_1 + \lambda_2 a_2 \mid \lambda_1, \lambda_2 \in \mathbb{Z}, |\lambda_1|+ |\lambda_2| \leq s \}.$$  It will be helpful to have the following facts available.

\begin{lem}  \label{lemma on index set}
For a given positive integer $s$, let $$I(s)=\{(\lambda_1,\lambda_2) \mid \lambda_1,\lambda_2 \in \mathbb{Z}, |\lambda_1|+|\lambda_2| \leq s \}.$$

\begin{enumerate}
  \item For a given integer $h$ with $|h| \leq s$, the number of ordered pairs $(\lambda_1, \lambda_2) \in I(s)$ with $\lambda_1 +\lambda_2 =h$ equals $s+1$ when $s$ and $h$ have the same parity, and $s$ otherwise.
  \item There are $(s+1)^2$ ordered pairs $(\lambda_1, \lambda_2) \in I(s)$ where $\lambda_1 +\lambda_2$ and $s$ have the same parity, and $s^2$ where $\lambda_1 +\lambda_2$ and $s$ have opposite parity. 

  \item There are $2s^2+2s+1$ elements in $I(s)$.
\end{enumerate}
\end{lem}

{\em Proof.}  By symmetry, we may assume that $h \geq 0$.  Observe that the `leftmost' integer point $(\lambda_1, \lambda_2)$ on the line $x+y=h$ within $|\lambda_1| + |\lambda_2| \leq s$ is 
$$\left( - \left \lfloor(s-h)/2 \right \rfloor , \; h + \left\lfloor (s-h)/2 \right \rfloor \right),$$
and the `rightmost' such point is 
$$\left( \left\lfloor (s+h)/2\right \rfloor ,  \; h- \left\lfloor (s+h)/2 \right \rfloor \right),$$ a set of $s+1$ points when $s$ and $h$ have the same parity, and $s$ points otherwise.  This establishes the first claim, from which our other two claims follow easily.  {\hfill $\Box$}

\begin{prop}

For all positive integers $s$ and $n$ with $$2 \leq n \leq 2s^2+2s+1,$$ the set $\{s,s+1\}$ is an $s$-spanning set in $\mathbb{Z}_n$.

\end{prop}

{\em Proof.}  The fact that $n$ cannot be larger than $2s^2+2s+1$ follows from Proposition \ref{diameter a(m,s)}.
We verify that for every $n$ in the range given above, we have $\langle \{s,s+1\} \rangle_s =\mathbb{Z}_n.$ 

Considering the elements of $$S=\{ \lambda_1 s + \lambda_2 (s+1) \mid \lambda_1, \lambda_2 \in \mathbb{Z}, |\lambda_1|+ |\lambda_2| \leq s \}$$ in $\mathbb{Z}$ (rather than $\mathbb{Z}_n$), we see that they
 lie in the interval $[-(s^2+s),(s^2+s)]$.  Since the index set $$I(s)=\{(\lambda_1,\lambda_2) \mid \lambda_1,\lambda_2 \in \mathbb{Z}, |\lambda_1|+|\lambda_2| \leq s \}$$ contains exactly $2s^2+2s+1$ elements (see Lemma \ref{lemma on index set}), it suffices to prove that no integer in $[-(s^2+s),(s^2+s)]$ can be written as an element of $S$ in two different ways.  

For that, suppose that $$\lambda_1s+\lambda_2(s+1)=\lambda_1's+\lambda_2'(s+1)$$ for some $(\lambda_1,\lambda_2) \in I(s)$ and $(\lambda_1',\lambda_2') \in I(s)$; without loss of generality, we can assume that $\lambda_2 \geq \lambda_2'$.  Our equation implies that $\lambda_2 - \lambda_2'$ is divisible by $s$ and is, therefore, equal to 0, $s$, or $2s$.  

If $\lambda_2-\lambda_2'=2s,$ then $\lambda_2=s$ and $\lambda_2'=-s$; since $|\lambda_1|+|\lambda_2| \leq s$ and $|\lambda_1'|+|\lambda_2'| \leq s$, this can only happen if $\lambda_1=\lambda_1'=0$, but this case leads to a contradiction with our equation.  

Assume, next, that $\lambda_2-\lambda_2'=s.$  Our equation then yields $\lambda_1-\lambda_1'=-s-1.$  In this case, we have $\lambda_1 \leq 0$, $\lambda_1' \geq 0$, $\lambda_2 \geq 0$, and $\lambda_2' \leq 0$, and we see that
$$|\lambda_1|+|\lambda_2|=-\lambda_1+\lambda_2=(s+1-\lambda_1')+(s+\lambda_2')=2s+1-(|\lambda_1'|+|\lambda_2'|) \geq s+1,$$ but that is a contradiction.

This leaves us with the case that $\lambda_2=\lambda_2'$, which implies $\lambda_1=\lambda_1'$, as claimed, and therefore the set $\{s,s+1\}$ is $s$-spanning in $\mathbb{Z}_n$.
{\hfill $\Box$}

\begin{prop}

For every positive integer $s$, the set $\{(0,1),(1,1)\}$ is an $s$-spanning set in $\mathbb{Z}_s \times \mathbb{Z}_{2s}$.

\end{prop}

{\em Proof.}  We need to prove that for any $(a,b) \in \mathbb{Z}_s \times \mathbb{Z}_{2s}$, we have $(a,b) \in \langle \{(0,1),(1,1)\} \rangle_s.$ 
Note that by symmetry we may assume that $0 \leq b \leq s$.  We consider three subcases.

When $b \geq a$, we let $\lambda_1=b-a$ and $\lambda_2=a$.  We then have 
$$\lambda_1 (0,1) + \lambda_2 (1,1) = (a,b)$$ and 
$$|\lambda_1|+|\lambda_2| = \lambda_1+\lambda_2=b \leq s,$$ so $(a,b) \in \langle \{(0,1),(1,1)\} \rangle_s.$

Next, when $2a-s \le b < a$, we again let $\lambda_1=b-a$ and $\lambda_2=a$.  In this case, we have
$$|\lambda_1|+|\lambda_2| = -\lambda_1+\lambda_2=2a-b \leq s,$$ so again $(a,b) \in \langle \{(0,1),(1,1)\} \rangle_s.$ 

Finally, when $b < 2a-s$, we let $\lambda_1=b-a+s$ and $\lambda_2=a-s$.  Here we have
$$\lambda_1 (0,1) + \lambda_2 (1,1) = (a-s,b),$$ where $a-s \equiv a$ mod $s$.  Furthermore,
$$|\lambda_1|+|\lambda_2| = \lambda_1-\lambda_2=(b -a+s)+(s-a) < s,$$ so again $(a,b) \in \langle \{(0,1),(1,1)\} \rangle_s.$
This completes our proof.  
{\hfill $\Box$}

\begin{thm}  \label{thm floor ceiling s-spanning}

Let $c$ and $s$ be positive integers so that the remainder of $s$ mod $c$ is at least $(c-1)/2$.  Let $u=\lfloor s/c \rfloor $, $v=\lceil s/c \rceil$, and $k=2 uv $.  Then $A=\left \{ (1,u) , (1, v) \right\}$ is an $s$-spanning set in $G=\mathbb{Z}_{c} \times \mathbb{Z}_{ck}$.

\end{thm}

{\em Proof.}  Let $H=\langle (1,u) \rangle$ be the subgroup of $G$ generated by the element $(1,u)$.  Since $H$ has order $2cv$, it has $cu$ cosets in $G$.  We claim that these cosets are $\lambda_2 \cdot (1,v) +H$, where $\lambda_2 = 0, 1, \ldots, cu-1$.  In order to prove this, we need to show that these $cu$ cosets are pairwise distinct.

Suppose, then, that $\mu_1 (1,v)$ and $\mu_2 (1,v)$ are in the same coset of $H$ for some $0 \leq \mu_1 , \mu_2 \le cu-1$.  This implies that $(\mu_2-\mu_1) \cdot (1,v)$ is an element of $H$, so there exists an integer $\mu$ for which $$(\mu_2-\mu_1) \cdot (1,v) = \mu \cdot (1,u).$$  From the first coordinates, this means that $\mu_2-\mu_1-\mu$ is divisible by $c$, and from the second coordinates we see that 
$(\mu_2-\mu_1)v - \mu u$, which, since $v=u+1$, equals $(\mu_2-\mu_1-\mu)u +(\mu_2-\mu_1),$ is divisible by $ck$ and thus by $cu$.  Combining these, we get that $\mu_2-\mu_1$ must also be divisible by $cu$, but that is only possible if $\mu_1=\mu_2$.

Let $g \in G$ be arbitrary.  According to our claim, we have integers $\lambda_1 \in [-cv+1,cv]$ and $\lambda_2 \in [0,cu-1]$ for which 
$$g=\lambda_1 (1,u) + \lambda_2 (1,v).$$
By assumption, we have integers $q$ (quotient) and $r$ (remainder) with $(c-1)/2 \le r \le c-1$ for which $s=qc+r$.  We then have
\begin{eqnarray} \label{eq with 2s} |\lambda_1|+|\lambda_2| \le cv + (cu-1) = c (q+1)+(cq-1) =2cq +c-1 = 2s +(c-2r-1) \le 2s. \end{eqnarray}  
If $|\lambda_1|+|\lambda_2| \le s$, then $ g \in \langle A \rangle_s$, and we are done.

Assume then that $|\lambda_1|+|\lambda_2| \ge s+1$. In this case, let $\lambda_1'=\lambda_1-\epsilon cv$ and $\lambda_2'=\lambda_2- cu$, where $\epsilon$ equals $1$ when $\lambda_1 \ge 0$, and $-1$ when $\lambda_1 < 0$.
We can compute that
$$\lambda_1' (1,u) + \lambda_2'(1,v)=(\lambda_1-\epsilon cv) (1,u) + (\lambda_2- cu) (1,v)=\lambda_1 (1,u) + \lambda_2 (1,v) - (\epsilon cv +cu, (\epsilon+1)cuv ),$$
and since $\epsilon +1$ equals $0$ or $2$, in $G$ we have $$\lambda_1' (1,u) + \lambda_2'(1,v)= \lambda_1 (1,u) + \lambda_2 (1,v) = g.$$

Finally, we show that $|\lambda_1'|+|\lambda_2'| \le s$.
Indeed, we have   
$$|\lambda_1'|+|\lambda_2'| = |\lambda_1-\epsilon cv| + |\lambda_2- cu| = cv - |\lambda_1| + cu - |\lambda_2|.$$
From (\ref{eq with 2s}) above we see that $cv+cu \le 2s+1$, and by our current assumption we have $|\lambda_1|+|\lambda_2| \ge s+1$.  Therefore,  
$|\lambda_1'|+|\lambda_2'| \le s$, as claimed, and our proof is now complete.
{\hfill $\Box$}

\begin{cor}  \label{corollary for min noncyclic}

Let $s$ be a positive integer for which $2s+1$ is composite, and let $p$ be any prime divisor of $2s+1$.  Then the largest noncyclic group with an $s$-spanning set of size $2$ has order at least $ 2s^2 +2s - (p^2-1)/2.$

\end{cor}

{\em Proof.}  Let us write $2s+1=pt$ for some integer $t \ge 2$.  Note that $p$ and $t$ must be odd, and we have $$s=\frac{t-1}{2}p+\frac{p-1}{2},$$ $\lfloor s/p \rfloor =(t-1)/2$, and $\lceil s/p \rceil=(t+1)/2$.  Let $G=\mathbb{Z}_p \times \mathbb{Z}_{pk}$ where $k=(t^2-1)/2$.  Then $G$ has order $ 2s^2 +2s - (p^2-1)/2$ and, according to Theorem \ref{thm floor ceiling s-spanning}, it has an $s$-spanning set of size $2$.
{\hfill $\Box$}

\section{$s$-regular groups and upper bounds}  \label{section -regular groups and upper bounds}

Recall that we say that the group $\mathbb{Z}_c \times \mathbb{Z}_{kc}$ is {\em $s$-regular} for some positive integer $s$, if it possesses an $s$-spanning set of the form $\{(1,u),(1,v)\}$ for some $u, v \in \mathbb{Z}_{kc}$.  In this section we provide some upper bounds for the largest possible order of an $s$-regular group.

\begin{thm}  \label{s regular prime c}

Let $G=\mathbb{Z}_{c} \times \mathbb{Z}_{ck}$ for some integers $c \geq 2$ and $k \geq 1$ that has an $s$-spanning set for some positive integer $s$. If $c$ is an odd prime, then $G$ is $s$-regular.

\end{thm}

{\em Proof.}  
First, we prove that $G$ has an $s$-spanning set of the form $A=\{(a,u),(b,v)\}$ for some $a,b \in \mathbb{Z}_{c}$ and $u, v \in \mathbb{Z}_{ck}$ with $u \not \equiv v $ mod $c$.  Note that if $\{(a,u),(b,v)\}$ is an $s$-spanning set, then so is $\{(a,u),(-b,-v)\}$, so if our claim were false, then we would have $u \equiv v \equiv-v$ mod $c$ for each $s$-spanning set $A=\{(a,u),(b,v)\}$ of $G$.  But since $c$ is an odd prime, that would imply that $u \equiv v \equiv0$ mod $c$, and that is a contradiction with $A$ being a spanning set.

Suppose then that $A=\{(a,u),(b,v)\}$ is an $s$-spanning in $G$ with $u \not \equiv v $ mod $c$.  We will prove that, in this case, $B=\{(1,u),(1,v)\}$ is also an $s$-spanning in $G$.  

Let $$I(s)=\{(\lambda_1,\lambda_2) \mid \lambda_1,\lambda_2 \in \mathbb{Z}, |\lambda_1|+|\lambda_2| \leq s \},$$ and consider the function $f: I(s) \rightarrow G$, given by $$(\lambda_1,\lambda_2) \mapsto \lambda_1 (a,u) + \lambda_2 (b,v).$$  Since $A$ is an $s$-spanning set, $f$ is surjective; we can thus find a subset $I_0$ of $I(s)$ on which the restriction $f_0 : I_0 \rightarrow G$ of $f$ is bijective.  We will prove that $g_0 : I_0 \rightarrow G$, given by $$(\lambda_1,\lambda_2) \mapsto \lambda_1 (1,u) + \lambda_2 (1,v)$$ is also bijective.  It suffices to show that $g_0$ is injective.  In order to do so, let us assume that the equation
$$\lambda_1 (1,u) + \lambda_2 (1,v) = \mu_1 (1,u) + \mu_2 (1,v)$$ holds in $G$
for some $(\lambda_1,\lambda_2)$ and $(\mu_1,\mu_2)$ of $I_0$.  We then have
\begin{eqnarray} \label{first eq mod c} 
\lambda_1  + \lambda_2 \equiv \mu_1  + \mu_2 \; \text{mod} \; c,
\end{eqnarray}
and 
\begin{eqnarray} \label{second eq mod ck}\lambda_1 u + \lambda_2 v \equiv \mu_1 u + \mu_2 v \; \text{mod} \; ck.
\end{eqnarray}
From (\ref{second eq mod ck}) we find that 
$$(\lambda_1 + \lambda_2)u+ \lambda_2(v-u) \equiv (\mu_1  + \mu_2)u + \mu_2 (v-u) \; \text{mod} \; ck;$$ considering (\ref{first eq mod c}), this gives us that
$$\lambda_2(v-u) \equiv  \mu_2 (v-u)  \; \text{mod} \; c.$$  Therefore, since $c$ is prime and  $u \not \equiv v $ mod $c$, we get that $\lambda_2 \equiv  \mu_2 $ mod $c.$  From (\ref{first eq mod c}) then we get $\lambda_1 \equiv  \mu_1 $ mod $c$, and thus 
\begin{eqnarray} \label{third eq mod c}\lambda_1 a + \lambda_2 b \equiv \mu_1 a + \mu_2 b \; \text{mod} \; c.
\end{eqnarray}
But $f_0$ is an injection, so  (\ref{third eq mod c}) and (\ref{second eq mod ck}) can hold together only when $(\lambda_1,\lambda_2)=(\mu_1,\mu_2)$.  This completes our proof.
{\hfill $\Box$}

\begin{thm}  \label{thm on s-regular max order by type}

Let $G=\mathbb{Z}_{c} \times \mathbb{Z}_{ck}$ for some integers $c \geq 2$ and $k \geq 1$, and suppose that $G$ is $s$-regular for some positive integer $s$. 
\begin{enumerate}
  \item If $c$ is even, then $|G| \leq 2s^2,$ and equality may occur only when $s \equiv 0$ mod $c$.
  \item If $c$ is odd and 
  \begin{itemize} 
    \item $s \equiv 0$ mod $c$, then $|G| \leq 2s^2+s$;
    \item $s \equiv (c-1)/2$ mod $c$, then $|G| \leq 2s^2 + 2s -(c^2-1)/2;$ 
    \item $c=3$ and $s \equiv 2$ mod $3$, then $|G| \leq 2s^2 +1$.
    \end{itemize}  In all other cases, $|G| < 2s^2$.
\end{enumerate}

\end{thm}

{\em Proof.}  Let $A=\{(1,u),(1,v)\}$ be an $s$-spanning set in $G$, in which case 
$$\{(\lambda_1 + \lambda_2, \lambda_1 u + \lambda_2 v) \mid (\lambda_1, \lambda_2) \in \mathbb{Z}^2, |\lambda_1| + |\lambda_2| \leq s\} =G.$$

Recall from Lemma \ref{lemma on index set} that, for a given integer $h$ with $|h| \leq s$, the number $\Lambda (h)$ of ordered pairs $(\lambda_1, \lambda_2) \in \mathbb{Z}^2$ with $|\lambda_1| + |\lambda_2| \leq s$ and $\lambda_1 +\lambda_2 =h$ equals $s+1$ when $s$ and $h$ have the same parity, and $s$ otherwise.   For each $i=0,1, \ldots, c-1$, we let $N(i)$ denote the sum of all $\Lambda (h)$  with $|h| \leq s$ and $h \equiv i$ mod $c$.    
Since $A$ is an $s$-spanning set in $G$, we must have
$$\{(\lambda_1 + \lambda_2, \lambda_1 u + \lambda_2 v) \mid (\lambda_1, \lambda_2) \in \mathbb{Z}^2, |\lambda_1| + |\lambda_2| \leq s\} = \mathbb{Z}_c \times \mathbb{Z}_{ck},$$ and thus for each $i=0,1, \ldots, c-1$,
$$\{(\lambda_1 + \lambda_2, \lambda_1 u + \lambda_2 v) \mid (\lambda_1, \lambda_2) \in \mathbb{Z}^2, |\lambda_1| + |\lambda_2| \leq s, \lambda_1 + \lambda_2 \equiv i \; \mbox{mod} \; c\} = \{i\} \times \mathbb{Z}_{ck}.$$
In particular,
\begin{eqnarray}\label{lambda i}
|G| \leq c N(i)
\end{eqnarray} must hold for each $i=0,1, \ldots, c-1$.
We now inspect (\ref{lambda i}) in several cases.  

When $c$ is even, we consider $i=r+1$, where $r$ is the nonnegative remainder of $s$ mod $c$ (with $r+1$ understood to be $0$ when $r=c-1$).  To compute $N(r+1)$, we need to identify the values of $h$ for which $|h| \leq s$ and $h \equiv r+1$ mod $c$, and these values are $h=s+1-jc$ for $j=1,2, \ldots, j_{\mathrm{max}}$, where
$$j_{\mathrm{max}}=\lfloor (2s+1)/c \rfloor = (2s-2r)/c + \lfloor (2r+1)/c \rfloor.$$  Note that $\lfloor (2r+1)/c \rfloor$ equals $0$ when $0 \leq r \leq c/2-1$, and $1$ when $c/2 \leq r \leq c-1$.  Since $s$ and $h=s+1-jc$ have different parity, and thus $\Lambda(h)=s$ for each $j$, we have
$$cN(r+1) = \left (2s-2r + c \lfloor (2r+1)/c \rfloor \right) s,$$ and we see that this quantity is at most $2s^2$, with equality when $r=0$ or $r=c/2$.  But $|G|$ must be divisible by $c^2$, so $|G|=2s^2$ may only occur when $r=0$.

When $c$ is odd and $r \leq c-2$, we consider $i=r+2$ (where $r+2$ is understood to be $0$ when $r=c-2$).  In this case the values of $h$ with $|h| \leq s$ and $h \equiv r+2$ mod $c$ are $h=s+2-tc$ for $t=1,2, \ldots, t_{\mathrm{max}}$, where
$$t_{\mathrm{max}}=\lfloor (2s+2)/c \rfloor = (2s-2r)/c + \lfloor (2r+2)/c \rfloor.$$  Note that $\lfloor (2r+2)/c \rfloor$ equals $0$ when $0 \leq r \leq (c-3)/2$, and $1$ when $(c-1)/2 \leq r \leq c-2$.  Since $s$ and $h=s+2-tc$ have the same parity when $t$ is even and different parity when $t$ is odd, we have
$$cN(r+2) = \left (s-r + c \lfloor (2r+2)/c \rfloor \right) s + (s-r)(s+1).$$
We consider four subcases.  For $r=0$, this quantity equals $2s^2+s$, and when $1 \leq r \leq (c-3)/2$, it equals $(s-r)s+(s-r)(s+1)$, which is less than $2s^2$.  When $(c+1)/2 \leq r \leq c-2$, we have 
$$cN(r+2) = \left (s-r + c \right) s + (s-r)(s+1) = (s-r)(2s+1)+cs \leq \left (s-(c+1)/2\right) (2s+1) + cs,$$ which again is less than $2s^2$.
Our last subcase is when $r=(c-1)/2$, when we have 
$$cN(r+2) = \left(s-(c-1)/2\right) (2s+1) + cs = 2s^2+2s-(c-1)/2.$$ 
Observe that when $s \equiv (c-1)/2$ mod $c$, then $2s^2+2s-(c^2-1)/2 \equiv 0$ mod $c^2$, and since $|G|$ must be divisible by $c^2$, this is the largest it can then be.

Our final case is when $c$ is odd and $r=c-1$; we consider $i=0$.  In this case the values of $h$ with $|h| \leq s$ and $h \equiv 0$ mod $c$ are $h=s+1-zc$ for $z=1,2, \ldots, z_{\mathrm{max}}$, where
$$z_{\mathrm{max}}=\lfloor (2s+1)/c \rfloor = (2s+2-c)/c.$$  Note that $z_{\mathrm{max}}$ is odd.  Since $s$ and $h=s+1-zc$ have the same parity when $z$ is odd and different parity when $z$ is even, we have 
$$N(0)=(z_{\mathrm{max}}-1)/2 \cdot s + (z_{\mathrm{max}}+1)/2 \cdot (s+1),$$ and thus
$$cN(0) = \left(s+1 -c \right)s + \left(s+1 \right)(s+1)=2s^2-(c-3)s+1,$$
which is less than $2s^2$, unless $c=3$.  This completes our proof.
{\hfill $\Box$}

\begin{thm}  \label{thm on s-regular max order}

For $s \geq 3$, every $s$-regular abelian group of rank 2 has order at most $2s^2+2s-4$, with equality if, and only if, it is isomorphic to $\mathbb{Z}_3 \times \mathbb{Z}_{(2s^2+2s-4)/3}$ for some $s \equiv 1$ mod $3$.  

\end{thm}

{\em Proof.}  Suppose that $G=\mathbb{Z}_c \times \mathbb{Z}_{ck}$ is an $s$-regular group of order at least $2s^2+2s-4$ for some integers $s \ge 3$, $c \geq 2$, and $k \geq 1$.  Then, by Theorem \ref{thm on s-regular max order by type}, $c$ is odd; furthermore,
\begin{itemize}
  \item $s \equiv 0$  mod $c$, and $2s^2+2s-4 \leq c^2k \leq 2s^2+s$, or
  \item $s \equiv (c-1)/2$  mod $c$, and $2s^2+2s-4 \leq c^2k \leq 2s^2+2s-(c^2-1)/2.$
\end{itemize}
Since the first possibility leads to no solutions, our claim is established.
{\hfill $\Box$}

\section{The groups $\mathbb{Z}_2 \times \mathbb{Z}_{2k}$}  \label{section groups 2 by 2k}

In this section we provide some results for $s$-spanning sets in groups of the form $\mathbb{Z}_2 \times \mathbb{Z}_{2k}$, including the exact value of the largest possible order of such groups.  Our general conjecture is as follows.

\begin{conj}

Let $k$ and $s$ be positive integers.  The group $\mathbb{Z}_2 \times \mathbb{Z}_{2k}$ has an $s$-spanning set of size $2$ if, and only if, one of the following holds:

\begin{itemize}
  \item $s$ is odd and $k \leq (s^2-1)/2$;
  \item $s$ is even and $k \leq (s^2-s)/2$;
  \item $s\equiv 0$ mod $4$, $k$ is even, and $(s^2-s)/2 < k \leq s^2/2$; or
  \item $s\equiv 2$ mod $4$, $k\equiv 2$ mod $4$, and $(s^2-s)/2 < k \leq s^2/2$.
\end{itemize} 
\end{conj}

\begin{thm}  \label{c=2 upper}

Let $k$ and $s$ be positive integers.  If the group $G=\mathbb{Z}_2 \times \mathbb{Z}_{2k}$ has an $s$-spanning set of size 2, then it has order at most $2s^2$.
\end{thm}

{\em Proof.}  If $G$ is $s$-regular, that is, it has an $s$-spanning set of the form $\{(1,u),(1,v)\}$ for some $u,v \in \mathbb{Z}_{2k}$, then our claim follows from Theorem \ref{thm on s-regular max order by type}.  Since no set of the form $\{(0,u),(0,v)\}$ can be a spanning set, we can suppose that $G$ has an $s$-spanning set of the form $A=\{(0,u),(1,v)\}$ for some $u,v \in \mathbb{Z}_{2k}$.  Therefore, each element in $\langle A \rangle_s$ has the form $(\lambda_2,\lambda_1 u + \lambda_2 v)$ for some integers $\lambda_1$ and $\lambda_2$ with $|\lambda_1|+|\lambda_2| \leq s$.   

Clearly, if $u$ and $v$ are both even, then $A$ can only generate elements of $G$ with an even second component, and thus it cannot be a spanning set.

Now if $u$ and $v$ are both odd, then $\lambda_1 u + \lambda_2 v$ has the same parity as $\lambda_1  + \lambda_2 $ does.  According to Lemma \ref{lemma on index set}, the parity of $\lambda_1  + \lambda_2 $ differs from the parity of $s$ exactly $s^2$ times (and agrees with it $(s+1)^2$ times) within $|\lambda_1|+|\lambda_2| \leq s$.  This means that $A$ cannot generate more than $s^2$ elements in $\mathbb{Z}_2 \times \mathbb{Z}_{2k}$ whose second coordinate has different parity than $s$ does, so if it is an $s$-spanning set, then we must have $|G| \leq 2s^2$.

Next, if $u$ is even and $v$ is odd, then for each element of $G$ that $A$ generates, its two components will have the same parity (namely, the parity of $\lambda_2$).  Therefore, $A$ cannot be a spanning set.

Finally, we prove that if $A=\{(0,u),(1,v)\}$ is an $s$-spanning set of $G$ with $u$ odd and $v$ even, then $A'=\{(1,u),(1,v)\}$ is an $s$-spanning set as well, meaning that $G$ is $s$-regular and thus has order at most $2s^2$ by Theorem \ref{thm on s-regular max order by type}.  

Let $a \in \mathbb{Z}_2$ and $b \in \mathbb{Z}_{2k}$.  Since $\langle A \rangle_s=G$, there are integers $\lambda_1$ and $\lambda_2$ with $|\lambda_1|+|\lambda_2| \leq s$ for which
$$ \lambda_1 (0,u) + \lambda_2 (1,v)=(\lambda_2,\lambda_1 u + \lambda_2 v) = (a',b),$$ where $a'=a$ if $b$ is even and $a'=1+a$ when $b$ is odd.  Observe that when $b$ is even, then (since $u$ is odd and $v$ is even), $\lambda_1$ must be even, and thus 
$$ \lambda_1 (1,u) + \lambda_2 (1,v)=(\lambda_1+\lambda_2,\lambda_1 u + \lambda_2 v) = (\lambda_2,\lambda_1 u + \lambda_2 v) = (a',b)=(a,b).$$ Similarly, when $b$ is odd, $\lambda_1$ must be odd, and 
$$ \lambda_1 (1,u) + \lambda_2 (1,v)= (1+a',b)=(a,b).$$
Therefore, $\langle A' \rangle_s=G$, as claimed, which completes our proof.
{\hfill $\Box$}

\begin{thm}  \label{thm s odd span upper}

Let $k$ and $s$ be positive integers with $s$ odd.  Then the group $G=\mathbb{Z}_2 \times \mathbb{Z}_{2k}$ has an $s$-spanning set of size 2 if, and only if, $k \leq (s^2-1)/2$.
\end{thm}

{\em Proof.}  The fact that $k \leq (s^2-1)/2$ is necessary follows from Theorem \ref{c=2 upper}.  We may also recall that Theorem \ref{thm floor ceiling s-spanning},  applied with $c=2$ and odd $s$, yields that with $u=(s-1)/2$ and $v=(s+1)/2$, the set $\{(1,u),(1,v)\}$ is an $s$-spanning set in $\mathbb{Z}_2 \times \mathbb{Z}_{s^2-1}$. This takes care of the case $k=(s^2-1)/2$.  

Our goal below is to show that the same values of $u=(s-1)/2$ and $v=(s+1)/2$ as {\em integers} provide an $s$-span for the set $\mathbb{Z}_2 \times \{0,1,2, \ldots, k\}$ when $k \leq (s^2-3)/2$.  (Here we consider the operation in the first component in $\mathbb{Z}_2$ and in the second component in $\mathbb{Z}$.)  By symmetry, we can then conclude that $\mathbb{Z}_2 \times \mathbb{Z}_{2k}$ has an $s$-spanning set of size 2 for $k \leq (s^2-3)/2$ as well.

We thus need to establish that for every $a \in \mathbb{Z}_2$ and $b \in \{0,1,2,\ldots,k\}$,
there are integers $\mu_1$ and $\mu_2$ with $|\mu_1|+|\mu_2| \leq s$ for which
\begin{eqnarray} \label{eqn with the mu coeffs}
\mu_1 (1,(s-1)/2) + \mu_2 (1,(s+1)/2)=(a,b)
\end{eqnarray} holds in $\mathbb{Z}_2 \times \{0,1,2, \ldots, k\}$.

Let us recall from the proof of Theorem \ref{thm floor ceiling s-spanning} (applied again for odd $s$ and $c=2$) that there are integers $\lambda_1 \in [-s,s+1]$ and $\lambda_2 \in [0,s-2]$ for which one of the following two sets of conditions holds:
\begin{enumerate}[(i)]
  \item $|\lambda_1|+|\lambda_2| \leq s$, and equation (\ref{eqn with the mu coeffs}) is satisfied in $\mathbb{Z}_2 \times \mathbb{Z}_{s^2-1}$ with $(\mu_1,\mu_2)=(\lambda_1,\lambda_2)$; or
  \item $|\lambda_1|+|\lambda_2| \geq s+1$, and equation (\ref{eqn with the mu coeffs}) is satisfied in $\mathbb{Z}_2 \times \mathbb{Z}_{s^2-1}$ with $(\mu_1,\mu_2)=(\lambda_1 - \epsilon (s+1),\lambda_2-s+1)$, where $\epsilon$ equals $1$ when $\lambda_1 \ge 0$, and $-1$ when $\lambda_1 < 0$.
\end{enumerate}

To complete our proof, we will show that, considering $$B=\mu_1 (s-1)/2 + \mu_2 (s+1)/2 $$ as an integer (rather than as an element of $\mathbb{Z}_{s^2-1}$), we have $B=b$.  Since $0 \leq b \leq k < (s^2-1)/2$, it suffices to prove that 
\begin{eqnarray} \label{B bounds} -(s^2-1)/2 \leq B < s^2-1. \end{eqnarray}

In case (i), we have $\mu_1=\lambda_1 \in [-s,s+1]$ and $\mu_2=\lambda_2 \in [0,s-2]$, from which (\ref{B bounds}) readily follows.  For case (ii), we find that
$$|\mu_1|+|\mu_2| = (s+1-|\lambda_1|) + (s-1 - |\lambda_2|) \leq 2s - (s+1) =s-1$$  and thus
$$|B| = |\mu_1 (s-1)/2 + \mu_2 (s+1)/2| \le |\mu_1|(s-1)/2 + |\mu_2| (s+1)/2 \le (|\mu_1|+|\mu_2|) (s+1)/2 \leq (s^2-1)/2,$$ so (\ref{B bounds}) holds again.  This completes our proof.
{\hfill $\Box$}

\begin{thm}  \label{thm s even span upper}

Let $k$ and $s$ be positive integers with $s$ even.  Then the group $G=\mathbb{Z}_2 \times \mathbb{Z}_{2k}$ has an $s$-spanning set of size 2 if $k \leq (s^2-s)/2$.
\end{thm}

{\em Proof.}  Like in the proof of Theorem \ref{thm s odd span upper}, we will prove that there are integers $u$ and $v$ for which the $s$-span of $A=\{(1,u),(1,v)\}$ contains the set $\mathbb{Z}_2 \times \{0,1,2, \ldots, k\}$.  (Here too we consider the operation in the first component in $\mathbb{Z}_2$ and in the second component in $\mathbb{Z}$.)  

Since $s-1$ is odd, from the proof of Theorem \ref{thm s odd span upper} we see that with $u=(s-2)/2$ and $v=s/2$, $\langle A \rangle_{s-1}$ contains the set 
$\mathbb{Z}_2 \times \{0,1,2,\ldots,k\}$ when $k \leq ((s-1)^2-1)/2$.  Therefore, $(1,v)+\langle A \rangle_{s-1}$ contains 
$\mathbb{Z}_2 \times \{s/2, s/2+1,\ldots,s/2+k\}$ for each $k \leq (s^2-2s)/2$, and since $(1,v) + \langle A \rangle_{s-1} \subseteq \langle A \rangle_{s}$, our claim follows for all $k \leq (s^2-s)/2$.
{\hfill $\Box$}


\begin{thebibliography}{99}


\bibitem{Baj:2018a} B. Bajnok, Additive Combinatorics: A Menu of Research Problems.  {\em CRC Press, Boca Raton}, 2018, xix+390 pp. 

\bibitem{BajBerJus:2022a} B. Bajnok, C. Berson, H. A. Just, On perfect bases in finite Abelian groups.  {\em Involve} {\bf 15} (2022), no. 3, 525--536.


\bibitem{Del:1895a} H. Delannoy, Emploi de l'\'echiquier pour la r\'esolution de certains probl\`emes de probabilit\'es.  {\em Assoc. Franc. Bordeaux} {\bf 24} (1895), 70--90.


\bibitem{Fio:2013a} M. A. Fiol, Comments on ``Extremal Cayley digraphs of finite abelian groups'' [Intercon. Networks 12 (2011), no. 1--2, 125--135].  {\em J. Inter. Net.} {\bf 14} (2013), no. 4, 1350016.
 
\bibitem{FioYebAleVal:1987a} M. A. Fiol, J. L. A. Yebra, I. Alegre, and M. Valero, A discrete optimization problem in local networks and data alignment.  {\em IEEE Trans. Comput.} {\bf 36} (1987), no. 6, 702--713.

 
\bibitem{GraSlo:1980a} R. L. Graham and N. J. A. Sloane, On additive bases and harmonious graphs.  {\em SIAM J. Algebraic Discrete Methods} {\bf 1} (1980), no. 4, 382--404.

\bibitem{Haa:2004a} H. Haanp\"a\"a, Minimum sum and difference covers of abelian groups.  {\em J. Integer Seq.} {\bf 7} (2004), no. 2, Article 04.2.6, 10 pp.
 

\bibitem{HsuJia:1994a} D. F. Hsu and X. Jia, Extremal problems in the construction of distributed loop networks.  \emph{SIAM J. Discrete Math.} {\bf 7} (1994), no. 1, 57--71.

\bibitem{MasSchJia:2011a} A.G. Mask, J. Schneider, X. Jia, Extremal Cayley digraphs of finite Abelian groups,
{\em J. Intercon. Networks} {\bf 12} (2011), no. 1--2, 125--135.

\bibitem{MorFioFab:1985a} P. Morillo, M. A. Fiol, and J. F\`abrega, The diameter of directed graphs associated to plane tessellations.  Tenth British combinatorial conference (Glasgow, 1985). {\em Ars Combin.} {\bf 20-A} (1985), A, 17--27.

\end{thebibliography}
\end{document}